\input amstex
\magnification=1200
\documentstyle{amsppt}
\loadeusm
\loadeufb
\loadbold
\def\qed{\hfill$\ssize\square$}
\predefine\accute{\'}
\redefine\'{\kern.05em{}}
\def\<{\langle}
\def\>{\rangle}

\def\im{\operatorname{im}}
\def\tr{\operatorname{tr}}

\def\D{\Cal D}

\def\H{\Cal H}
\def\V{\Cal V}

\def\d{\delta}
\def\k{\kappa}
\def\s{\sigma}

\def\lie#1{L\lower.8ex\hbox{$\ssize#1$}\'}
\def\nab#1{\nabla\kern-.2em\lower.8ex\hbox{$\ssize#1$}\'}
\def\diffop#1#2{#1\kern-.1em\lower.8ex\hbox{$\ssize#2$}\'}
\def\nabsq#1#2{\nabla^2\kern-.3em\lower.8ex\hbox{$\ssize#1,#2$}\,} 
\def\diffopsq#1#2#3{#1^2\kern-.3em\lower.8ex\hbox{$\ssize#2,#3$}\,}
\def\midstar{\'\hbox{$*$}\'}
\def\uperp#1{#1\raise1.1ex\hbox{$\sssize\bot$}}
\PSAMSFonts

\NoBlackBoxes
\nologo
\TagsOnRight
\rightheadtext{Harmonic Contact Metric Structures, and Submersions}
\pagewidth{16truecm}
\pageheight{24.3truecm}
\hoffset=-0.5truecm
\voffset=-0.5truecm
\linespacing{1.5}
\topmatter
\title
Harmonic Contact Metric Structures, and Submersions
\endtitle
\author 
E\. Vergara-Diaz
\&\
C\. M\. Wood
\endauthor
\address 
School of Mathematics, Trinity College Dublin, Dublin 2, Ireland.
\newline\indent
Department of Mathematics, University of York, Heslington, York
Y010 5DD, UK.
\endaddress
\email 
evd\@maths.tcd.ie,
cmw4\@york.ac.uk 
\endemail
\dedicatory
In memoriam, Professor James Eells, 1926--2007.
\enddedicatory
\abstract
We study harmonic almost contact structures in the context of contact metric manifolds, and an analysis is carried out when such a manifold fibres over an almost Hermitian manifold, as exemplified by the Boothby-Wang fibration.  Two types of almost contact metric warped products are also studied, relating their harmonicity to that of the almost Hermitian structure on the base or fibre.
\endabstract
\endtopmatter
\document
\head
1. Introduction
\endhead
In recent years the study of {\sl harmonic unit vector fields\/} on Riemannian manifolds has attracted considerable interest: see \cite{6} for a fairly 
contemporary survey.  Such fields are named for the fact that they are harmonic sections of the unit tangent bundle, and as such may be considered ``optimal'' when compared with nearby unit vector fields.  Geometric interest is added by the fact that their total bending is stationary \cite{11}.  Although there is no general existence theory for harmonic vector fields, many examples arise as characteristic, or Reeb, fields on contact or almost contact manifolds.  In fact contact metric manifolds with harmonic Reeb field $\xi$ have been dubbed {\sl $H$-contact\/} in \cite{9}, and characterized as those where $\xi$ is an eigenvector of the Ricci operator.  They include {\sl $K$-contact\/} manifolds \cite{2}. However, when comparing almost contact metric structures one would also like to take into account the geometry of the {\sl hyperplane bundle\/} $\D\to M$ orthogonal to the Reeb field, in particular its induced Hermitian structure.  For this reason, in \cite{10} we introduced the idea of a {\sl harmonic almost contact metric structure\/} on an orientable Riemannian manifold $(M^{2n+1},g)$, defined as follows. Let $N\to M$ be the fibre bundle with fibre $SO(2n+1)/U(n)$ associated to the principal bundle of $g$-orthonormal tangent frames of $M$ viz\. the odd-dimensional analogue of the twistor bundle in Hermitian geometry.  Then an almost contact metric structure on $M$ is parametrized by a unique section $\s$ of $N$, and is said to be {\sl harmonic\/} if $\s$ is a harmonic section.  In \cite{10} we analysed the geometry of the homogeneous bundle $N\to M$ and thereby showed that $\s$ is a harmonic section precisely when the following {\it two\/} equations are satisfied:
$$
\gather
\tau(\xi)+\tfrac12\'J\'T(\phi)=0,
\tag1.1 \\
\tau(J)=0.
\tag1.2
\endgather
$$
Here, $\phi$ is the fundamental $(1,1)$ tensor of the almost contact structure, satisfying:
$$
\phi^2=-I+\eta\otimes\xi,
$$
where $\eta$ is the $1$-form dual to $\xi$, and $J$ is the almost complex structure in $\D=\im\phi$
\linebreak
obtained by restricting $\phi$.  If $\bar\nabla$ is the connection in $\D$ obtained by orthogonally projecting the Levi-Civita connection $\nabla$ of $(M,g)$, then $\tau(\xi)$ and $\tau(J)$ are defined as follows:
$$
\tau(\xi)=\nabla^*\nabla\xi-|\nabla\xi|^2\xi,
\qquad
\tau(J)=[\bar\nabla^*\bar\nabla J,J],
$$
where $\nabla^*\nabla$ and $\bar\nabla^*\bar\nabla$ are the {\sl rough Laplacians\/} of $(TM,\nabla)$ and $(\D,\bar\nabla)$, respectively:
$$
\nabla^*\nabla\xi=-\tr\nabla^2\xi=-\nabsq{E_i}{E_i}\xi,
\qquad
\bar\nabla^*\bar\nabla J
=-\tr\bar\nabla^2J=-\bar\nabsq{E_i}{E_i}J,
$$
for any local orthonormal frame $(E_i)$ of $TM$.  As will be the case throughout the paper, the summation signs in these formul\ae\ have been omitted.  Notice that $\tau(\xi)$ is the $\D$-component of $\nabla^*\nabla\xi$.  Finally, $T(\phi)$ is the following section of $\D$:
$$
T(\phi)=\tr(\bar\nabla J\otimes\nabla\xi)
=\bar\nab{E_i}J(\nab{E_i}\xi).
$$
Equation (1.2) is simply the condition for $J$ to be a harmonic almost complex structure in $\D$, in the sense of \cite{13}.  On the other hand, $\tau(\xi)=0$ is the equation for $\xi$ to be a harmonic vector field \cite{15,\,16}, so (1.1) is perhaps more subtle than one might expect.  In particular, an almost contact metric structure with harmonic Reeb field is harmonic only if the additional condition  $T(\phi)=0$ holds.  In this paper we refer to equations (1.1) and (1.2) as the {\sl first\/} and {\sl second harmonic section equations,} respectively.  
\par
A general analytic existence theory for equations (1.1) and (1.2) is currently far from being within reach; the best available existence results for harmonic sections only apply to bundles with compact negatively curved fibres \cite{12,\,14}.  We therefore seek geometric techniques for constructing solutions.  In \cite{10} we studied the case where $M^{2n+1}$ is an isometrically immersed hypersurface of an almost Hermitian manifold $\tilde M^{2n+2}$, and the almost contact structure is induced by the ambient almost complex structure.  The relationship between equations (1.1) and (1.2) and the almost Hermitian geometry of $\tilde M$ is in general rather complicated, but becomes more tractable when $\tilde M$ is a K\"ahler or nearly K\"ahler manifold, or $\tilde M=M\times\Bbb R$.  In this paper we investigate the `dual' problem, when $M^{2n+1}$ fibres over a $2n$-dimensional almost Hermitian manifold $\hat M^{2n}$.  In the first part of the paper we confine attention to the case when $M$ is a contact metric manifold.  We begin \S2 with a characterization (Proposition 2.1) of the first harmonic section equation in terms of:
$$
h=\tfrac12\lie\xi\phi,
\tag1.3
$$
which is an important structural tensor in contact metric geometry \cite{3}.  
We then note an analogy between contact metric manifolds and $(1,2)$-symplectic manifolds in almost Hermitian geometry (Lemma 2.1), which guides us to a characterization of the second harmonic section equation in terms of the {\sl \midstar Ricci curvature\/} $\bar\rho^*$ of $\D$, which is defined:
$$
\bar\rho^*(X,Y)=g(\bar R(X,F_i)JF_i,JY),
\quad\text{for all $X,Y\in\D$,}
\tag1.4
$$
where $\bar R$ denotes the curvature of the vector bundle $(\D,\bar\nabla)$:
$$
\bar R(X,Y)=\bar\nabsq XY-\bar\nabsq YX,
$$
and $(F_i)$ is a local orthonormal frame of $\D$.  We show (Proposition 2.2) that the second harmonic section equation is equivalent to $J$-invariance (or equivalently, symmetry) of $\bar\rho^*$.  The analogous result in almost Hermitian geometry is that a $(1,2)$-symplectic structure is harmonic if and only if its \midstar Ricci curvature is symmetric \cite{13}.  Combining Propositions 2.1 and 2.2 yields our main result, Theorem 2.1, which is a characterization of harmonic contact metric structures.  In particular, an $H$-contact structure is harmonic if and only if $h$ is co-closed (when viewed as a $TM$-valued $1$-form) and $\bar\rho^*$ is symmetric, and a $K$-contact structure is harmonic precisely when $\bar\rho^*$ is symmetric.  For practical purposes it is helpful to recast this in terms of the {\sl \midstar Ricci curvature of $M$,} defined:
$$
\rho^*(X,Y)=g(R(X,E_i)\phi E_i,\phi Y),
\quad\text{for all $X,Y\in TM$,}
\tag1.5
$$
where $R$ is the Riemann curvature tensor of $(M,g)$.  We then show (Theorem 2.2) that an $H$-contact structure is harmonic precisely when $\rho^*$ is symmetric, or equivalently $\phi$-invariant.  As an example, we show (Theorem 2.3) that all contact metric structures satisfying the $(\k,\mu)$-nullity condition of \cite{4} are harmonic.  Such structures include the unit tangent bundles of spaces of constant curvature.
\par
In \S3 we take the analogy of \S2 one step further, to the situation where there exists a Riemannian submersion $\pi\colon M\to\hat M$ onto an almost Hermitian manifold $\hat M^{2n}$, which intertwines the contact metric structure of $M$ with the almost complex structure of $\hat M$.  Then $M$ is necessarily a $K$-contact manifold and $\hat M$ is almost K\"ahler, and $\pi$ exists under certain topological conditions, with the classical construction of the Boothy-Wang fibration \cite{5}.  We use Theorem 2.1 to prove that the contact metric structure on $M$ is harmonic precisely when the almost Hermitian structure on $\hat M$ is harmonic (Theorem 3.1).  We then consider the same setup when $M$ is merely an almost contact manifold, in the special case where $M=\hat M\times_f\Bbb R$, the warped product equipped with the induced almost contact structure.  The result here (Theorem 3.2) is that the almost contact structure is harmonic if and only if the almost Hermitian structure is harmonic, provided that the gradient vector $\nabla f$ is K\"ahler null.
Finally we consider the reversed warped product $M=\Bbb R\times_f\check M$, with the almost contact structure induced by an almost Hermitian structure on $\check M^{2n}$, or more generally where $M$ is locally of this form; for example if $M$ is a Kenmotsu manifold, which was studied in \cite{10}.  We show (Theorem 3.3) that the almost contact structure is harmonic if and only if the almost Hermitian structure is harmonic, provided the latter is cosymplectic, or the warping function is constant.  It is interesting (Propositions 3.1 and 3.3) that in both warped products the Reeb field is harmonic, for any warping function $f$, which we believe constitutes a new family of examples of harmonic unit vector fields.
\par
\bigskip
\head
2. Harmonic contact metric structures
\endhead
There are a number of useful relations in the tensor algebra of a contact metric manifold, documented in \cite{3}.  Firstly:
$$
\nab\xi\phi=0,
\tag2.1
$$ 
from which it follows that $\xi$ is geodesic: 
$$
\nab\xi\xi=0.
\tag2.2
$$
Also, the tensor $h$ defined in (1.3) is symmetric, trace-free, anti-commutes with $\phi$, and verifies:
$$
\nab X\xi=-\phi X-\phi\'hX.
\tag2.3
$$
We may regard $h$ as a $TM$-valued $1$-form on $M$, and form its co-derivative:
$$
\d h=-\nab{E_i}h(E_i).
$$
Since $h\xi=0$, and $h$ is symmetric, $h$ is $\D$-valued.  Our first result shows that $\d h$ is in fact a section of $\D$, and appears naturally in the first harmonic section equation.  During the proof, and hereafter, we will often denote the Riemannian metric by $\<\,,\>$ instead of $g$. 
\proclaim{Proposition 2.1}
On a contact metric manifold, there is the following identity:
$$
\d h=\phi\,\nabla^*\nabla\xi-T(\phi).
\tag2.4
$$
The first harmonic section equation is equivalent to:
$$
\tau(\xi)=\phi\'\d h.
\tag2.5
$$
The first harmonic section equation is verified on an $H$-contact manifold if and only if $h$ is co-closed; in particular, it is always verified for a $K$-contact structure.
\endproclaim
\demo{Proof} 
First we recall from \cite{10} the relation:
$$
\bar\nabla J=\nabla\phi-\<\nabla\phi,\xi \>\xi.
\tag2.6
$$
We also note from (2.2) that:
$$
hX=\phi\,\nab X\xi-X,
\quad\text{for all $X\in\D$.}
\tag2.7
$$
Now suppose for convenience that a local orthonormal frame $(F_i)$ of $\D$ has been constructed by $\bar\nabla$-parallel translation of an orthonormal basis of the fibre $\D_x$ of $\D$ over $x\in M$, along radial geodesics; thus $\bar\nabla F_i(x)=0$.  Then at $x$ we also have:
$$
\align
\nab{F_i}F_i
&=\<\nab{F_i}F_i,\xi\>\xi 
=-\<F_i,\nab{F_i}\xi\>\xi
=\<F_i,\phi F_i+\phi hF_i\>\xi,
\quad\text{by (2.2)} \\
&=-\<\phi,h\>\xi=0,
\endalign
$$
since $h$ (resp\. $\phi$) is symmetric (resp\. skew-symmetric).  Therefore:
$$
\allowdisplaybreaks
\align
-\d h(x)
&=\nab{F_i}h(F_i)+\nab\xi h(\xi)
=\nab{F_i}(hF_i),
\quad\text{since also $h\xi=0$ and $\xi$ is geodesic,} \\
&=\nab{F_i}(\phi\,\nab{F_i}\xi),
\quad\text{by (2.7)} \\
&=\nab{F_i}\phi(\nab{F_i}\xi)+\phi\,\nab{F_i}\nab{F_i}\xi \\
&=\bar\nab{F_i}J(\nab{F_i}\xi)
+\<\nab{F_i}\phi(\nab{F_i}\xi),\xi\>\xi
-\phi\,\nabla^*\nabla\xi,
\quad\text{by (2.6)} \\
&=T(\phi)-\phi\,\nabla^*\nabla\xi
+\<\nab{F_i}\xi,\phi\,\nab{F_i}\xi\>\xi,
\quad\text{by (2.2)} \\
&=T(\phi)-\phi\,\nabla^*\nabla\xi,
\quad\text{since $\phi$ is skew-symmetric.}
\endalign
$$
This establishes (2.4).  Applying $\phi$ to both sides yields:
$$
\phi\'\d h=-\tau(\xi)-J\'T(\phi),
$$ 
which may be used with (1.1) to obtain (2.5).  An $H$-contact structure verifies $\tau(\xi)=0$, in which case (2.5) is equivalent to $\d h=0$.  Finally, a $K$-contact structure is $H$-contact \cite{9}, and is characterized by $h=0$; hence (2.5) is verified.
\qed
\enddemo
We now turn to the second harmonic section equation.  Our approach is guided by the following key observation, that $\bar\nabla J$ is $J$-anti-invariant.  
\proclaim{Lemma 2.1}
In a contact metric manifold, the Hermitian vector bundle $(\D,\bar\nabla,J)$ satisfies:
$$
\bar\nab{JX}J(JY)=-\bar\nab XJ(Y),
\quad\text{for all $X,Y\in\D$.}
$$
\endproclaim
\demo{Proof} 
The K\"ahler form for $J$ is the restriction to $\D$ of the $2$-form $\Phi(X,Y)=g(X,\phi Y)$.  Suppose $X,Y,Z\in\D$.  Then:
$$
\align
d\Phi(X,Y,Z)&=\nab X\Phi(Y,Z)+\nab Y\Phi(Z,X)+\nab Z\Phi(X,Y) \\
&=\<Y,\nab X\phi(Z)\>+\<Z,\nab Y\phi(X)\>+\<X,\nab Z\phi (Y)\> \\
&=\<Y,\bar\nab XJ(Z)\>+\<Z,\bar\nab YJ(X)\>+\<X,\bar\nab ZJ(Y)\>,
\quad\text{by (2.6).}
\endalign
$$
One can then establish the following (remarkable) identity, which we leave as an exercise for the reader:
$$
\multlinegap{1cm}
\multline
d\Phi(X,Y,Z)-d\Phi(X,JY,JZ)+d\Phi(JX,JY,Z)+d\Phi(JX,Y,JZ) \\
=-2\<\bar\nab XJ(Y)+\bar\nab{JX}J(JY),Z\>.
\endmultline
$$
The result follows since $\Phi=d\eta$ for a contact metric structure; thus $d\Phi=0$.
\qed
\enddemo
Lemma 2.1 may also be derived from Lemma 7.3 of \cite{3}; however, Lemma 2.1 is {\it not\/} equivalent to the $\phi$-anti-invariance of $\nabla\phi$.
\par
It follows from (2.1) and (2.6) that in addition to the relation of Lemma 2.1 we have:
$$
\bar\nab\xi J=0.
\tag2.9
$$  
Contact metric manifolds may therefore be regarded as odd-dimensional analogues of $(1,2)$-symplectic manifolds.  It should be noted however that this analogy is not a characterization.  For example, the Kenmotsu almost contact metric structure on the warped product $M=\Bbb R\times_f\check M^{2n}$, where $f(t)=e^t$ and $\check M$ is a K\"ahler manifold, satisfies $\bar\nabla J=0$, but is not a contact structure since $\D=T\check M$ which is clearly integrable.  We will see other examples of this in \S3.
\par
In order to derive a characterization of the harmonicity of the almost complex structure $J$ analogous to that of \cite{13} for $(1,2)$-symplectic structures, we now prove a technical lemma, which introduces the curvature tensor.
\proclaim{Lemma 2.2}
Let $F$ be an element of $\D$, the contact subbundle of a contact metric manifold.  Then:
$$
\big[\bar\nabsq FFJ-2\bar R(F,JF)+\bar\nabsq{JF}{JF}J,J\big]\,
=\,4\,\bar\nab{\bar\nabla J(F,F)}J.
$$
\endproclaim
\demo{Proof} 
Suppose that $F,X\in\D_x$, with extensions to local sections of $\D$ such that
$\bar\nabla X(x)=0=\bar\nabla F(x)$.  Note first that in this case, since the
Levi-Civita connection is torsion free, we get using (2.3):
$$
\align
[F,JF]&=\nab F(JF)-\nab{JF}F \\
&=\bar\nab F(JF)+\<\nab F(JF),\xi\>\xi-\<\nab{JF}F,\xi\>\xi \\
&=\bar\nab FJ(F)-\<JF,\nab F\xi\>\xi+\<F,\nab{JF}\xi\>\xi \\
&=\bar\nab FJ(F)+\<JF,JF+JhF\>\xi-\<F,J^2F+JhJF\>\xi \\
&=\bar\nabla J(F,F)+2|F|^2\xi,
\quad\text{since $h$ anticommutes with $\phi$.}
\tag2.10
\endalign
$$
Also, by Lemma 2.1:
$$
\bar\nab{JF}JF
=\bar\nab{JF}J(F)
=\bar\nab FJ(JF)
=-J\,\bar\nabla J(F,F).
\tag2.11
$$
We then note that:
$$
\align
\bar\nabsq FFJ
&=\bar\nab F\bar\nab FJ-\bar\nab{\nab FF}J
=\bar\nab F\bar\nab FJ-\<\nab FF,\xi\>\,\bar\nab\xi J \\
&=\bar\nab F\bar\nab FJ,
\quad\text{by (2.9),}
\tag2.12
\endalign
$$
and similarly, using (2.11) and Lemma 2.1:
$$
\bar\nabsq{JF}{JF}J
=\bar\nab{JF}\bar\nab{JF}J+J\,\bar\nab{\bar\nabla J(F,F)}J.
\tag2.13
$$
\indent
We now begin the main calculation.  Using the Leibnitz rule and (2.12), we have:
$$
\bar\nabsq FFJ(JX)
=\bar\nab F\bar\nab FJ(JX)
=\bar\nab F(\bar\nab FJ(JX))-\diffopsq{(\bar\nabla J)}FF(X),
$$
and the first term on the right hand side may be expanded using the curvature tensor as follows:
$$
\align
\bar\nab F(\bar\nab FJ(JX)) 
&=\bar\nab F(\bar\nab{JF}J(X)),
\quad\text{by Lemma 2.1 } \\
&=\bar\nab F\bar\nab{JF}(JX)-\bar\nab F(J\,\bar\nab{JF}X) \\
&=\bar\nab F\bar\nab{JF}(JX)-J\,\bar\nab F\bar\nab{JF}X \\
&=\bar\nab{JF}\bar\nab F(JX)+\bar\nab{[F,JF]}(JX)+\bar R(F,JF)(JX) \\
&\qquad-J\,\bar\nab{JF}\bar\nab FX-J\bar R(F,JF)X \\
&=\bar\nab{JF}\bar\nab F(JX)-J\,\bar\nab{JF}\bar\nab FX 
+\bar\nab{\bar\nabla J(F,F)}J(X)+[\bar R(F,JF),J]X,
\endalign
$$
using (2.10) and (2.9). The calculation continues:
$$
\align
\bar\nab F(\bar\nab FJ(JX))
&-\bar\nab{\bar\nabla J(F,F)}J(X)-[\bar R(F,JF),J]X \\
&=\bar\nab{JF}(\bar\nab FJ(X))+\bar\nab{JF}(J\,\bar\nab FX)
-J\,\bar\nab{JF}\bar\nab FX \\
&=-\bar\nab{JF}(\bar\nab{JF}J(JX)),\quad\text{by Lemma 2.1 again,} \\
&=-\bar\nab{JF}\bar\nab{JF}J(JX)-\bar\nab{JF}J(\bar\nab{JF}(JX)) \\
&=-\bar\nabsq{JF}{JF}J(JX)+\bar\nab{\bar\nabla J(F,F)}J(X)
-\diffopsq{(\bar\nabla J)}{JF}{JF}(X),
\quad\text{by (2.13).}
\endalign
$$
It follows from Lemma 2.1 once again that:
$$
\diffopsq{(\bar\nabla J)}{JF}{JF}(X)
=\bar\nab FJ\circ J\circ\bar\nab FJ(JX)
=\diffopsq{(\bar\nabla J)}FF(X).
$$
Therefore:
$$
\multlinegap{1cm}
\multline
\bar\nabsq FFJ (JX)+\bar\nabsq{JF}{JF}J(JX) \\
=[\bar R(F,JF),J]X
+2\,\bar\nab{\bar\nabla J(F,F)}J(X)
-2\,\diffopsq{(\bar\nabla J)}FF(X).
\endmultline
\tag2.14
$$  
Since $(\bar\nabla J)^2$ commutes with $J$, whereas $[\bar R(F,JF),J]$ anticommutes with $J$, replacing $X$ by $JX$ in (2.14) and applying $J$ to the resultant equation yields:
$$
\multlinegap{1cm}
\multline
J\,\bar\nabsq FFJ (X)+J\,\bar\nabsq{JF}{JF}J(X) \\
=-[\bar R(F,JF),J]X
-2\,\bar\nab{\bar\nabla J(F,F)}J(X)
-2\,\diffopsq{(\bar\nabla J)}FF(X).
\endmultline
\tag2.15
$$ 
Subtraction of (2.15) from (2.14) yields the result.
\qed
\enddemo
For any local orthonormal frame $(F_i)$ of the hyperplane bundle $\D$, we define analogously to \cite{13}:  
$$
\bar\delta J=-\tr_\D\bar\nabla J=-\bar\nab{F_i}J(F_i).
$$ 
Then it is an immediate consequence of Lemma 2.1 that $\bar\delta J=0$ for a contact metric structure.  Furthermore, on any almost contact metric manifold, the curvature tensor $\bar R$ of $\D$ is related to the $\D$-component $R_\D$ of the curvature tensor $R$ of $(M,g)$ by the following identity \cite{10}:
$$
\bar R(X,Y)Z=R_\D(X,Y)Z+r(\nab X\xi,\nab Y\xi)Z,
\tag2.16
$$
for all $X,Y,Z\in\D$, where $r$ is the following curvature-type tensor:
$$ 
r(u,v)w=\<v,w\>u-\<u,w\>v.
\tag2.17
$$
\proclaim{Proposition 2.2}
A contact metric structure verifies the second harmonic section equation if and only if the \midstar Ricci curvature of $\D$ is symmetric (or, equivalently, $J$-invariant). 
\endproclaim
\demo{Proof}
We first note the identity:
$$
[\bar\nabla^*\bar\nabla J,J]=-[\bar R(F_i,JF_i),J],
\tag2.18
$$
which follows from Lemma 2.2, using $\bar\delta J=0$ and $\bar\nabsq\xi\xi J=0$,
where the latter is a consequence of (2.9) and (2.2).  Now suppose $Z,W\in\D$.  Then, using Bianchi's first identity, on summation:
$$
R(F_i,JF_i)Z=-2R(Z,F_i)JF_i.
\tag2.19
$$
Therefore by (2.16) and (2.3), and making extensive use of the symmetry of $h$ and anticommutativity of $h$ and $J$:
$$
\align
\<\bar R(F_i,JF_i)Z,JW\> 
&=-2\<R(Z,F_i)JF_i,JW\>-\<r(JF_i+JhF_i,F_i+hF_i)Z,JW\> \\
&=-2\<R(Z,F_i)JF_i,JW\>-2\<Z,W\>-2\<hZ,hW\> \\
&=-2\<\bar R(Z,F_i)JF_i,JW\>+2\<r(Z+hZ,F_i+hF_i)F_i,W\> \\
&\qquad
-2\<Z,W\>-2\<hZ,hW\> \\
&=-2\'\bar\rho^*(Z,W)+4(n-1)\<Z+hZ,W\>-4\<hZ,hW\>,
\endalign
$$
where for the final equation we have also used the facts that $h$ is trace-free and $h\xi=0$.  Therefore, since $h$ is symmetric:
$$
\<[\bar R(F_i,JF_i),J](JZ),JW\> 
=2\'\bar\rho^*(Z,W)-2\'\bar\rho^*(W,Z),
$$
and the result follows from (2.18) and (1.2).
\qed
\enddemo
It is interesting to note that the characterization of harmonic $(1,2)$-symplectic structures in \cite{13} was obtained in a completely different way to Proposition 2.2, using a technique of Lichnerowicz.  Combining Propositions 2.1 and 2.1 yields:
\proclaim{Theorem 2.1} 
\flushpar
\rom{(1)}\quad
A contact metric structure is harmonic if and only if $\tau(\xi)=\phi\'\d h$ and the \midstar Ricci curvature of $\D$ is symmetric.
\flushpar
\rom{(2)}\quad
An $H$-contact structure is harmonic if and only if $h$ is co-closed and
the \midstar Ricci curvature of $\D$ is symmetric.
\flushpar
\rom{(3)}\quad
A $K$-contact structure is harmonic if and only if the \midstar Ricci curvature of $\D$ is symmetric.
\endproclaim
From the proof of Proposition 2.2 it is possible to extract the equation:
$$
\align
\bar\rho^*(Z,W)
&=\rho^*(Z,W)+(2n-1)\<Z,W\> \\
&\qquad
+2(n-1)\<hZ,W\>-2\<hZ,hW\>,
\tag2.20
\endalign
$$
where the \midstar Ricci curvature $\rho^*$ of $M$ is defined in (1.5).  Therefore the symmetry of $\bar\rho^*$ is equivalent to the symmetry of $\rho^*(\D,\D)$.  However, whereas $\rho^*(X,\xi)=0$ for all $X\in TM$, we have:
\proclaim{Lemma 2.3}
On a contact metric manifold:
$$
\rho^*(\xi,Z)=-\<\d h,JZ\>,
\quad\text{for all $Z\in\D$.}
$$
\endproclaim
\demo{Proof}
We first recall the following curvature identity for contact metric manifolds \cite{3,\,Lemma 7.4}:
$$
\<R(\xi,X)Y,Z\>=-\nab X\Phi(Y,Z)-\<X,\nab Y(\phi h)Z\>+\<X,\nab Z(\phi h)Y\>,
$$
from which it follows that for all $Z\in\D$ we have:
$$
\rho^*(\xi,Z)=-\nab{F_i}\Phi(JF_i,JZ)-\<F_i,\nab{JF_i}(\phi h)JZ\>
+\<F_i,\nab{JZ}(\phi h)JF_i\>.
$$
Calculating each term on the right hand side in turn:
$$
\align
\nab{F_i}\Phi(JF_i,JZ)\,
&=\<JF_i,\nab{F_i}\phi(JZ)\>
=-\<\bar\nab{F_i}J(JF_i),JZ\>,
\quad\text{by (2.6)} \\
&=-\<\bar\d J,Z\>=0,
\quad\text{by Lemma 2.1.}
\endalign
$$
Similarly, using the anticommutativity of $h$ and $J$:
$$
\align
\<F_i,\nab{JF_i}(\phi h)JZ\>
&=\<F_i,\nab{JF_i}\phi(hJZ)+\phi\,\nab{JF_i}h(JZ)\> \\
&=-\<\bar\nab{F_i}J(F_i),hZ\>-\< \nab{JF_i}h(JF_i),JZ\>,
\quad\text{by (2.6) and Lemma 2.1,} \\
&=\<\bar\d J,hZ\>+\<\d h,JZ\>
=\<\d h,JZ\>,
\quad\text{by Lemma 2.1.}
\endalign
$$
To compute the final term, we assume that $(F_i)$ has been chosen as in the proof of Proposition 2.1, and then note that for all $X\in TM$ we have:
$$
\tr\nab Xh
=\<\nab Xh(F_i),F_i\>
=\<\nab X(hF_i),F_i\>
=X.\<hF_i,F_i\>
=X.\tr h=0,
\tag2.21
$$
since $h\xi=0$ and $h$ is trace-free.  Hence:
$$
\align
\<F_i,\nab{JZ}(\phi h)JF_i\>
&=-\<F_i,\nab{JZ}(h\phi)JF_i\>
=\<F_i,\nab{JZ}h(F_i)-h\,\nab{JZ}\phi(JF_i)\> \\
&=\<F_i,h\,\bar\nab ZJ(F_i)\>,
\quad\text{by (2.21), (2.6) and Lemma 2.1,} \\
&=\,\<h,\bar\nab ZJ\>=\<h,\nab Z\phi\>=0,
\endalign
$$
since $h$ (resp\. $\nab Z\phi$) is symmetric (resp\. skew-symmetric).
\qed
\enddemo
\proclaim{Theorem 2.2}
An $H$-contact structure on $M$ is harmonic if and only if the \midstar Ricci curvature of $M$ is symmetric.
\endproclaim
As an example, we consider the {\sl $(\k,\mu)$-manifolds\/} introduced in \cite{4}.  These are the contact metric manifolds whose curvature satisfies:
$$
R(X,Y)\xi=(\k+\mu\'h)r(X,Y)\xi,
\quad\text{for all $X,Y\in TM$,}
\tag2.22
$$
where $\k,\mu$ are constants, and $r$ is the curvature-type tensor defined in (2.17).  It was shown in \cite{9} that $(\k,\mu)$-manifolds are $H$-contact.  Furthermore in \cite{4} it was shown that (2.22) determines $R$ completely (for a contact metric manifold), which enables us to analyse the \midstar Ricci curvature.
\proclaim{Theorem 2.3}
A contact metric structure satisfying the $(\k,\mu)$-nullity condition is harmonic.
\endproclaim
\demo{Proof}
We note first that since $r(\D,\D)\xi=0$ it follows from (2.22) that $\rho^*(\xi,Z)=0$ for all $Z\in\D$; therefore $\rho^*(\xi,Z)=\rho^*(Z,\xi)$.  Now for all $Z,W\in\D$ we have the identity:
$$
2\'\rho^*(Z,W)-2\'\rho^*(W,Z)
=\<[R(F_i,JF_i),J](JZ),JW\>,
$$
which follows from (2.19).  The curvature identity \cite{4,\,Lemma 3.2} may be recast in the following succinct way:
$$
\multline
\<[R(X,Y),\phi]Z,W\>
=\<[(1+h)r(X,Y)(1+h),\phi]Z,W\> \\
+(1-\k)\<r(X,Y)\xi,\phi\'r(Z,W)\xi\>
+(1-\mu)\<h\'r(X,Y)\xi,\phi\'r(Z,W)\xi\>,
\endmultline
$$
from which it follows that:
$$
\<[R(F_i,JF_i),J](JZ),JW\> 
=\<[(1+h)r(F_i,JF_i)(1+h),J](JZ),JW\>.
$$
Now, using the anticommutativity of $J$ and $h$:
$$
\multlinegap{0pt}
\multline
[(1+h)r(F_i,JF_i)(1+h),J]
=(1+h)r(F_i,JF_i)(1+h)J\,
-\,J(1+h)r(F_i,JF_i)(1+h) \\
=(1+h)(1-h)\,
-\,(1+h)J(1+h)J\,
+\,J(1+h)J(1+h)\,
+\,J(1+h)J(1+h) \\
=2(1+h)(1-h)\,
-\,2(1-h)(1+h)=0.
\endmultline
$$
Therefore $\rho^*$ is symmetric, and the result follows from Theorem 2.2.
\qed
\enddemo
\bigskip
\head
3. Submersive almost contact structures
\endhead
Initially let $M^{2n+1}$ be an almost contact metric manifold.  We say that the almost contact metric structure is {\sl submersive\/} if there exists an almost Hermitian manifold $(\hat M^{2n},\hat g,\hat J)$ and a Riemannian submersion $\pi\colon(M,g)\to(\hat M,\hat g)$ such that the almost complex structures in $\D$ and $T\hat M$ are compatible:
$$
d\pi(JZ)=\hat Jd\pi(Z),\quad\text{for all $Z\in\D$.}
\tag3.1
$$ 
We refer to $\hat J$ as the {\sl projected\/} almost Hermitian structure, and
henceforward, for simplicity, usually make no notational distinction between $J$ and $\hat J$, denoting the latter by $J$.  We also denote both $g$ and $\hat g$ by $\<\,,\>$.  Mixing the terminology of almost contact geometry and Riemannian submersions \cite{8}, tangent vectors to $M$ in the direction of $\xi$ are {\sl vertical,} whereas elements of $\D$ are {\sl horizontal.}  Recall also that a vector field $X$ on $M$ is said to be {\sl basic\/} if and only if $X$ is horizontal and $\pi$-related to a vector field $\hat X$ on $\hat M$: $\pi_*X=\hat X$.  
We will denote by $\V$ and $\H$ the orthogonal projections of $TM$ onto the vertical and horizontal subbundles. Thus, for all $X\in TM$:
$$
\V X=\eta(X)\xi,
\quad\text{and}\quad
\H X=X_{\D}=\bar{X},
\quad\text{where}\quad 
X=\bar X+\eta(X)\xi.
$$
We will utilize O'Neill's structure tensor $A$, defined:
$$
A_XY=\V(\nab{\H X}\H Y)+\H(\nab{\H X}\V Y),
\quad\text{for all $X,Y\in TM$.}
$$ 
Finally we recall Lemma 2 and Theorem 2 of \cite{8}, which in our context become:
\proclaim{Lemma 3.1}
If $X$ and $Y$ are sections of $\D$, then 
$$
A_XY=\tfrac12\'\eta([X,Y])\xi.
$$
\endproclaim
\proclaim{Lemma 3.2}
If $X,Y,Z,H$ are elements of $\D$, then:
$$
\align
g(R(X,Y)Z,H)
&=\hat g(\hat R(\hat X,\hat Y)\hat Z,\hat H)-2\,g(A_XY,A_ZH) \\
&\qquad+g(A_YZ,A_XH)+g(A_ZX,A_YH).
\endalign
$$
\endproclaim
\par
Now, if $M$ is a compact regular $K$-contact manifold, then $M$ is submersive, via the Boothby-Wang fibration \cite{5}.  Conversely, suppose that $M$ is a submersive contact metric manifold.  Suppose $X\in\D$ is extended to a local basic field.  Then $\phi X$ is also basic, by (3.1), and therefore:
$$
\align
2\pi_*(hX)
&=\pi_*\'\lie\xi\phi(X)
=\pi_*[\xi,\phi X]-\pi_*\phi[\xi,X] \\
&=[\pi_*\xi,\pi_*\phi X]-J\pi_*[\xi,X],
\quad\text{by (3.1),} \\
&=-J[\pi_*\xi,\pi_*X]=0.
\endalign
$$
Since $h$ is $\D$-valued (on a contact metric manifold) it follows that $h=0$.
Thus $M$ is necessarily $K$-contact, and $\hat M$ is therefore almost K\"ahler.
\proclaim{Theorem 3.1}
A submersive contact metric structure is harmonic if and only if the projected almost Hermitian structure is harmonic.
\endproclaim
\demo{Proof}
First, since the contact metric manifold is necessarily $K$-contact ($h=0$), 
it follows from Lemma 3.1 that for all horizontal tangent vectors $Z,W$ we have:
$$
\align
2A_ZW
&=\<\nab XW-\nab WZ,\xi\>\xi=-\<W,\nab Z\xi\>\xi+\<X,\nab Z\xi\>\xi \\
&=\<W,JZ\>\xi-\<Z,JW\>\xi=2\<W,JZ\>\xi,
\quad\text{by (2.3).}
\endalign
$$
Therefore:
$$
A_ZW=\<JZ,W\>\xi.
\tag3.2
$$
Now, since $h=0$, it follows from (2.20) that for all $X,Y\in\D$ we have:
$$
\align
\bar\rho^*(X,Y)
&=\<R(X,F_i)JF_i,JY\>+(2n-1)\<X,Y\>,
\intertext{and then from Lemma 3.2:}
&=\<\hat R(\hat X,\widehat{F_i})J\widehat{F_i},J\hat Y\>-2\<A_XF_i,A_{JF_i}JY\> \\
&\qquad+\<A_{F_i}JF_i,A_XJY\>+\<A_{JF_i}X,A_{F_i}JY\>+(2n-1)\< X,Y\>,
\intertext{and finally using (3.2):}
&=\hat\rho^*(\hat X,\hat Y)-2\<JX,F_i\>\<J^2F_i,JY\> \\
&\qquad+\<JF_i,JF_i\>\<JX,JY\>+\<J^2F_i,X\>\<JF_i,JY\>+(2n-1)\<X,Y\> \\
&=\hat\rho^*(\hat X,\hat Y)+4n\<X,Y\>.
\endalign
$$
Therefore the \midstar Ricci curvature of $\D$ is symmetric if and only if the \midstar Ricci curvature $\hat\rho^*$ of $\hat M$ is symmetric, which by \cite{13} is equivalent to the almost complex structure of $\hat M$ being harmonic, since $\hat M$ is an almost K\"ahler manifold.  The result then follows from Theorem 2.1.
\qed
\enddemo
There are of course submersive almost contact metric structures which are not contact.  As an example, we consider the warped product $M=\hat M\times_f\Bbb R$, where $f\colon\hat M\to\Bbb R$ is a strictly positive function.  The induced almost contact metric structure is determined by the stipulation that $\xi=f^{-1}\partial_t$ (where for convenience we are identifying $f\circ\pi$ with $f$) and $\D=T\hat M$ equipped with the almost Hermitian structure of $\hat M$.  In order to analyse the harmonic section equations we recall the following elementary aspects of warped product geometry \cite{1,\,Lemma 7.3}.
\proclaim{Lemma 3.3}
Let $M$ be the warped product $\hat M\times_f\check M$.  Let $X,Y$ be projectable horizontal (ie\. basic) vector fields, and $V,W$ projectable vertical vector fields.  Then:
\flushpar
\rom{(1)}\quad
$\nab XV=\nab VX=\<X,f^{-1}\'\nabla f\>\'V$\rom;
\flushpar
\rom{(2)}\quad
$\H(\nab VW)=-\<V,W\>f^{-1}\'\nabla f$ and $\;\V(\nab VW)=\check\nab{V}W$\rom;
\flushpar
\rom{(3)}\quad
$\nab XY=\hat\nab{X}Y$.
\flushpar
Here $\<\,,\>$ denotes the warped metric, and projectable tangent vectors on $M$ are identified with their projections to $\hat M$ or $\check M$ as appropriate.  Furthermore $\nabla f$ denotes the horizontal lift of the gradient field of $f$, or equivalently the gradient of $f\circ\pi$.
\endproclaim
\proclaim{Proposition 3.1}
On the warped product almost contact manifold $\hat M\times_f\Bbb R$ the Reeb field is harmonic, and the first harmonic section equation is verified.
\endproclaim
\demo{Proof} 
First we show that the rough Laplacian of the Reeb field is: 
$$
\nabla^*\nabla\xi=f^{-2}\'|\nabla f|^2\'\xi,
\tag3.3
$$
so $\xi$ is a harmonic vector field.  To see this, suppose $X\in\D_x$ has been extended to a local basic vector field with $\nabla X(x)=0$, which is possible by Lemma 3.3.  Then by Lemma 3.3: 
$$
\nab X\xi
=(X.f^{-1})\partial_t+f^{-1}\'\nab X\partial_t
=-f^{-2}(X.f)\partial_t+f^{-1}(f^{-1}(X.f))\partial_t=0.
\tag3.4
$$
This implies that $\nabsq XX\xi=0$.  On the other hand, using that $f$ does not depend on $t$, and $\partial_t$ is a unit vector field on $\Bbb R$ with geodesic integral curves, it follows from Lemma 3.3 that: 
$$
\nab\xi\xi=f^{-2}\,\nab{\partial_t}\partial_t
=-f^{-3}\<\partial_t,\partial_t\>\nabla f
=-f^{-1}\'\nabla f,
\tag3.5
$$
and then:
$$
\align
\nabsq\xi\xi\xi 
&=\nab\xi\nab\xi\xi-\nab{\nab\xi\xi}\xi
=-\nab\xi(f^{-1}\'\nabla f)+f^{-1}\'\nab{\nabla f}(f^{-1}\partial_t) \\
&=f^{-2}\big(\nab{\nabla f}\partial_t-\nab{\partial_t}(\nabla f)
-|\nabla f|^2\xi\big)
=-f^{-2}\'|\nabla f|^2\'\xi,
\quad\text{by Lemma 3.3.} \\
\endalign
$$
This establishes (3.3), and hence $\tau(\xi)=0$.  We now compute:
$$
\phi\,\nab\xi\phi(\nab\xi\xi)
=-f^{-1}\'\phi\,\nab\xi\phi(\nabla f) 
=-f^{-1}\'\phi\big(\nab\xi(\phi\,\nabla f)-\phi\,\nab\xi(\nabla f)\big)=0,
$$
since $\nab\xi(\nabla f)$ and $\nab\xi(\phi\,\nabla f)$ are vertical, by Lemma 3.3.  It follows from (2.6) and (3.4) that $T(\phi)=0$.  Therefore (1.1) is verified.
\qed
\enddemo
Since we are no longer in a position to use the curvature results of \S2,
in order to analyse the second harmonic section equation we begin with the following general observation.
\proclaim{Lemma 3.4}
For any submersive almost contact metric structure, if $X,Y$ are elements of $\D$ then:
$$
\pi_*\bar\nab XJ(Y)=\hat\nab{\hat X}J(\hat Y).
$$
\endproclaim
\demo{Proof} 
We extend $Y$ to a basic vector field.  It follows from (3.1) that $JY$ is also basic, with $\widehat{JY}=J\hat Y$.  Furthermore, by \cite{8,\,Lemma 1}:
$$
\pi_*\'\nab XY=\hat\nab{\hat X}\hat Y.
\tag3.6
$$
Therefore:
$$
\align
\pi_*(\bar\nab XJ(Y)) 
&=\pi_*\bar\nab X(JY)-\pi_*\'J\,\bar\nab XY
=\hat\nab{\hat X}\widehat{JY}-J\pi_*\bar\nab XY,
\quad\text{by (3.1)} \\
&=\hat\nab{\hat X}(J\hat Y)-J\,\hat\nab{\hat X}\hat Y
=\hat\nab{\hat X}J(\hat Y).
\tag"\qed"
\endalign
$$
\enddemo
\proclaim{Proposition 3.2}
On the warped product almost contact manifold $\hat M\times_f\Bbb R$ the second harmonic section equation is verified if and only if the almost Hermitian structure on $\hat M$ satisfies:
$$
[\hat\nabla^*\hat\nabla J,J]+2f^{-1}J\,\hat\nab{\nabla f}J=0.
$$
\endproclaim
\demo{Proof}
We note first that if $Z\in\D_x$ is extended to a local basic field then since $\phi Z$ is also basic it follows from Lemma 3.3 that:
$$
\nab\xi\phi(Z)
=\nab\xi(\phi Z)-\phi\,\nab\xi Z
=\<\phi Z,f^{-1}\nabla f\>\xi.
$$
Therefore by (2.6):
$$
\bar\nab\xi J=0.
\tag3.7
$$
It then follows from (3.5) that:
$$
\bar\nabsq\xi\xi J
=-\bar\nab{\nab\xi\xi}J
=f^{-1}\'\bar\nab{\nabla f}J,
$$
and hence by Lemma 3.4:
$$
\pi_*\'\bar\nabsq\xi\xi J
=f^{-1}\'\hat\nab{\nabla f}J.
\tag3.8
$$
Now suppose that $X,Y\in\D_x$ are also extended to local basic fields.  It then follows from Lemma 3.4 that, for example, $\bar\nab YJ(Z)$ is also basic.  Let us further suppose that $Y,Z$ have been extended such that $\nabla Y(x)=0=\nabla Z(x)$ (cf\. the proof of Proposition 3.1).  Then at $x$:
$$
\align
\pi_*\'\bar\nabsq XYJ(Z)
&=\pi_*\'\bar\nab X\bar\nab YJ(Z)
=\pi_*\'\bar\nab X(\bar\nab YJ(Z)) \\
&=\hat\nab{\hat X}(\hat\nab{\hat Y}J(\hat Z)),
\quad\text{by (3.6) and Lemma 3.4,} \\
&=\hat\nabsq{\hat X}{\hat Y}J(\hat Z),
\quad\text{by Lemma 3.3.}
\tag3.9
\endalign
$$
Combining (3.8) and (3.9) yields:
$$
\pi_*[\bar\nabla^*\bar\nabla J,J]
=[\hat\nabla^*\hat\nabla J,J]+2f^{-1}J\,\hat\nab{\nabla f}J,
$$
and the result follows on comparison with (1.2).
\qed
\enddemo
It is perhaps worth noting that in contrast to Lemma 3.4, the attractive equation (3.9) does not hold for arbitrary submersive almost contact structures; it is in fact a consequence of $J$ being parallel along the Reeb field.
We also note that if $\hat M$ is $(1,2)$-symplectic then it follows from Lemma 3.4 and (3.7) that $\bar\nabla J$ has the same symmetries as those of a contact metric structure (Lemma 2.1 and (2.9)), although $M$ is no longer a contact manifold.
\par
Combining Propositions 3.1 and 3.2 gives our second main result on submersive almost contact metric structures.
\proclaim{Theorem 3.2}
The induced almost contact structure on $M=\hat M\times_f\Bbb R$ is harmonic if and only if the almost Hermitian structure on $\hat M$ satisfies:
$$
[\hat\nabla^*\hat\nabla J,J]+2f^{-1}J\,\hat\nab{\nabla f}J=0.
$$
If $\hat M$ is a K\"ahler manifold then the almost contact  structure on $M$ is harmonic for all warping functions $f$.  If the almost Hermitian structure on $\hat M$ is harmonic (for example, if $\hat M$ is a nearly K\"ahler manifold) then the almost contact structure on $M$ is harmonic if and only if $\nabla f$ is K\"ahler null.  
\endproclaim
It follows from Theorem 3.2 that if $\hat M$ is a {\sl strict\/} nearly K\"ahler manifold \cite{7} (for example, the six dimensional sphere) then the warped product almost contact structure on $\hat M\times_f\Bbb R$ is never harmonic.  We conclude by showing that this situation can be corrected if the warping is reversed, so that $M=\Bbb R\times_f\check M$ where $f\colon\Bbb R\to\Bbb R$ is strictly positive and $\check M^{2n}$ is an almost Hermitian manifold.  In this case the induced almost contact metric structure is specified by $\xi=\partial_t$, and $\D=T\check M$ equipped with the almost Hermitian structure of $\check M$.  Notice however that the almost contact structure is no longer submersive,
and $\xi$ is now horizontal, whereas $\D$ is vertical.  It follows from Lemma 3.3 that $\xi$ is geodesic.
\proclaim{Lemma 3.5}
The induced almost contact structure on $\Bbb R\times_f\check M$ verifies $\nab\xi\phi=0$, and:
$$
\bar\nab XJ(Y)=\check\nab XJ(Y),
\quad\text{for all $X,Y\in\D$.}
$$
\rom(All projections have been notationally omitted.\rom)
\endproclaim
\demo{Proof}
Since $\xi$ is geodesic, we have $\nabla\phi(\xi,\xi)=0$.  Furthermore if $Y\in\D$ is extended to a projectable vertical field then $\phi Y$ is also projectable vertical and so by Lemma 3.3:
$$
\nab\xi\phi(Y)
=\nab\xi(\phi Y)-\phi\,\nab\xi Y
=f^{-1}f'\phi Y-\phi(f^{-1}f'\'Y)=0.
$$
Now by (2.6) and Lemma 3.3, taking into account the identification of $\D$ and $T\check M$:
$$
\align
\bar\nab XJ(Y)
&=\V\,\nab X\phi(Y)
=\V\,\nab X(\phi Y)-\phi\V\,\nab XY \\
&=\check\nab X(JY)-J\,\check\nab XY
=\check\nab XJ(Y).
\tag"\qed"
\endalign
$$
\enddemo
For our next result, recall that an almost Hermitian manifold is said to be {\sl cosymplectic\/} if its K\"ahler form is co-closed.  The class of such manifolds includes all $(1,2)$-symplectic manifolds, and hence all almost K\"ahler and nearly K\"ahler manifolds.
\proclaim{Proposition 3.3}
On the warped product almost contact manifold $\Bbb R\times_f\check M$ the Reeb field is harmonic, and the first harmonic section equation is verified if and only if $\check M$ is a cosymplectic almost Hermitian manifold or $f$ is constant.
\endproclaim
\demo{Proof}
First we prove that the vector field $\xi$ is harmonic.  Since $\xi$ is geodesic, we have $\nabsq\xi\xi\xi=0$.  Now if $X\in\D_x$ is extended to a projectable vertical field with $\check\nabla X(x)=0$ then it follows from Lemma 3.3 that:
$$
\nab XX=-f^{-1}f'\'|X|^2\xi,
\quad\text{and}\quad
\nab X\xi=f^{-1}f'X,
\tag3.10
$$
noting that $\nabla f=f'\xi$.  Therefore, since $\xi$ is geodesic: 
$$
\align
\nabsq XX\xi 
&=\nab X\nab X\xi
=\nab X(f^{-1}f'X)= X.(f^{-1}f')\'X+f^{-1}f'\,\nab XX \\
&=-f^{-2}(f')^2|X|^2\xi.
\endalign
$$
It therefore follows that:
$$
\nabla^*\nabla \xi
=2nf^{-2}(f')^2\xi,
$$
and hence $\xi$ is harmonic.  Next, we prove that $T(\phi)$ is identically zero and the result will follow.  Let $(F_i)$ be a projectable vertical local orthonormal frame.  Then since $\xi$ is geodesic:
$$
\align
T(\phi)
&=\bar\nab{F_i}(\nab{F_i}\xi)
=f^{-1}f'\,\bar\nab{F_i}J(F_i),
\quad\text{by (3.10),} \\
&=f^{-1}f'\,\check\nab{F_i}J(F_i),
\quad\text{by Lemma 3.5,} \\
&=f^{-1}f'\'\check\d J=0,
\endalign
$$
if $f$ is constant or $\check M$ is cosymplectic.
\qed
\enddemo
\proclaim{Theorem 3.3}
Suppose that $\check M$ is a cosymplectic almost Hermitian manifold, or $f$ is constant.  Then the induced almost contact structure on $\Bbb R\times_f\check M$ is harmonic if and only if the almost Hermitian structure on $\check M$ is harmonic.  In particular, if $\check M$ is a nearly K\"ahler manifold then the induced almost contact structure is harmonic.
\endproclaim
\demo{Proof} 
We note first that by Lemma 3.5 and (2.6):
$$
\bar\nab\xi J=0.
\tag3.11
$$
Therefore since $\xi$ is geodesic:
$$
\bar\nabsq\xi\xi J
=\bar\nab\xi\bar\nab\xi J=0.
$$
Now suppose $X$ is as in the proof of Proposition 3.3.  Then:
$$
\align
\bar\nabsq XXJ
&=\bar\nab X\bar\nab XJ-\bar\nab{\nab XX}J \\
&=\bar\nab X\bar\nab XJ,
\quad\text{by (3.10) and (3.11),} \\
&=\V\,\nab X\check\nab XJ,
\quad\text{by Lemma 3.5,} \\
&=\check\nab X\check\nab XJ,
\quad\text{by Lemma 3.3,} \\
&=\check\nabsq XXJ,
\quad\text{since $\check\nabla X(x)=0$.}
\endalign
$$
It follows that:
$$
\bar\nabla^*\bar\nabla J=\check\nabla^*\check\nabla J,
$$
and hence (1.2) is verified precisely when the almost Hermitian structure on $\check M$ is harmonic.  Combining this with Proposition 3.3 yields the result.
\qed
\enddemo
Inspecting equation (3.11) and Lemma 3.5 shows that once again we are able to construct examples of almost contact structures where $\bar\nabla J$ has the same symmetries as those of a contact metric manifold.

\enddocument
\end